# NONSUBJECTIVE PRIORS VIA PREDICTIVE RELATIVE ENTROPY REGRET

By Trevor J. Sweeting,[1] Gauri S. Datta[2] and Malay Ghosh[3]

*University College London, University of Georgia and University of Florida*

We explore the construction of nonsubjective prior distributions in Bayesian statistics via a posterior predictive relative entropy regret criterion. We carry out a minimax analysis based on a derived asymptotic predictive loss function and show that this approach to prior construction has a number of attractive features. The approach here differs from previous work that uses either prior or posterior relative entropy regret in that we consider predictive performance in relation to alternative nondegenerate prior distributions. The theory is illustrated with an analysis of some specific examples.

**1. Introduction.** There is an extensive literature on the development of objective prior distributions based on information loss criteria. Bernardo [5] obtains reference priors by maximizing the Shannon mutual information between the parameter and the sample. These priors are maximin solutions under relative entropy loss; see, for example, [3, 8] for further analysis, discussion and references. In regular parametric families the reference prior for the full parameter is Jeffreys' prior. It is argued in [5], however, that when nuisance parameters are present, then the appropriate reference prior should depend on which parameter(s) are deemed to be of primary interest. This dependence on parameters of interest is mirrored in the approach to prior development via minimization of coverage probability bias; see, for example, [11, 23, 25] for further aspects of this approach.

In the present paper we explore the construction of nonsubjective prior distributions via predictive performance. It is possible to use Bernardo's ap-

Received March 2003; revised July 2005.
[1]Supported in part by EPSRC Grant GR/R24210/01.
[2]Supported in part by NSF Grants DMS-00-71642 and SES-02-41651 and NSA Grant MDA904-03-1-0016.
[3]Supported in part by NSF Grant SES-99-11485.

proach to obtain reference priors for prediction. However, as shown in [5], this program turns out to be equivalent to obtaining the reference prior for the full parameter, which produces Jeffreys' prior in regular problems. Further analysis along these lines is carried out in [17]. Datta et al. [12] explore prior construction using predictive probability matching, which is shown to produce sensible prior distributions in a number of standard examples. In the present article we follow Bernardo [5] and Barron [3] by taking an information-theoretic approach and using an entropy-based risk function. However, here we focus on the posterior predictive relative entropy regret, as opposed to the prior predictive relative entropy regret used by these authors. Our starting point is the predictive information criterion introduced by Aitchison [1], which was also discussed by Akaike [2] as a criterion for the selection of objective priors. We depart from these and other authors by taking a more Bayesian viewpoint, in that we are less concerned here with performance in repeated sampling but rather with performance in relation to alternative prior specifications. The main aim of the paper is to search for uniform, or impartial, minimax priors under an associated predictive loss function. These priors are also maximin, or least favorable, which can be interpreted here as giving rise to minimum information predictive distributions.

The organization of the paper is as follows. We start in Section 2 by defining the posterior predictive regret, which measures the regret when using a posterior predictive distribution under a particular prior in relation to the posterior predictive distribution under an alternative proper prior. We define a related predictive loss function and argue that this is a suitable criterion for the comparison of alternative prior specifications. We discuss informally the results in Section 6 on impartial, minimax and maximin priors under a large sample version of this loss function. We also give a definition of the predictive information in a prior distribution. Throughout we make connections with standard quantities that arise in information theory. In Section 3 we relate posterior predictive regret and loss to prior predictive regret and loss and in Section 4 we obtain the asymptotic behavior of the posterior predictive regret, which is obtained via an analysis of the higher-order asymptotic behavior of the prior predictive regret. The higher-order analysis carried out in Section 5, which is of independent interest, leads to expressions for the asymptotic forms of the posterior predictive regret, predictive information and predictive loss. In Section 6 we investigate impartial minimax priors under our asymptotic predictive loss function. It turns out that these priors also minimize the asymptotic information in the predictive distribution. In the case of a single real parameter, Jeffreys' prior turns out to be minimax. However, in dimensions greater than one, the minimax solution need not be Jeffreys' prior. The theory is illustrated with an analysis of some specific examples, and some concluding remarks are given in Section 7.

There are a number of appealing aspects of the proposed Bayesian predictive approach to prior determination. First, since the focus is on prediction, there is no need to specify a set of parameters deemed to be of interest. Second, difficulties associated with improper priors are avoided in the formulation of posterior predictive, as opposed to prior predictive, criteria. Third, the minimax priors identified in Section 6 arise as limits of proper priors. Fourth, these minimax priors are also maximin, or least favorable for prediction, which can be interpreted here as minimizing the predictive information contained in a prior. Finally, and importantly, the same asymptotic predictive loss criterion emerges regardless of whether one is considering prediction of a single future observation or a large number of future observations.

**2. Posterior predictive regret and impartial priors.** Consider a parametric model with density $p(\cdot|\theta)$ with respect to a $\sigma$-finite measure $\mu$, where $\theta = (\theta^1, \ldots, \theta^p)$ is an unknown parameter in an open set $\Theta \subset \mathcal{R}^p$, $p \geq 1$. Let $p^\pi(x) = \int p(x|\theta)\, d\pi(\theta)$ be the marginal density of $X$ under the prior distribution $\pi$ on $\Theta$, where both $\pi$ and $p^\pi$ may be improper. Let $\Pi$ be the class of prior distributions $\pi$ satisfying $p^\pi(X) < \infty$ a.s. $(\theta)$ for all $\theta \in \Theta$. That is, $\pi \in \Pi$ if and only if $P^\theta(\{X : p^\pi(X) < \infty\}) = 1$ for all $\theta \in \Theta$.

We suppose that $X$ represents data to be observed and $Y$ represents future observations to be predicted. Denote by $p^\pi(y|x)$ the posterior predictive density of $Y$ given $X = x$ under the prior $\pi \in \Pi$. Let $\Omega \subset \Pi$ be the class of all proper prior distributions on $\Theta$. For $\pi \in \Pi$ and $\tau \in \Omega$, define the *posterior predictive regret*

$$(2.1) \qquad d_{Y|X}(\tau, \pi) = \iint \log\left\{\frac{p^\tau(y|x)}{p^\pi(y|x)}\right\} p^\tau(x,y)\, d\mu(x)\, d\mu(y).$$

We note that $d_{Y|X}(\tau, \pi)$ is the conditional relative entropy, or expected Kullback–Leibler divergence, $D(p^\tau(Y|X)\|p^\pi(Y|X))$, between the predictive densities under $\pi$ and $\tau$. See, for example, the book by Cover and Thomas [10] for definitions and properties of the various information-theoretic quantities that arise in this work. It follows from standard results in information theory that the quantity $d_{Y|X}(\tau, \pi)$ always exists (possibly $+\infty$) and is nonnegative. It is zero when $\pi = \tau$ and is therefore the expected regret under the loss function $-\log p^\pi(y|x)$ associated with using the predictive density $p^\pi(y|x)$ when $X$ and $Y$ arise from $p^\tau(x)$ and $p^\tau(y|x)$, respectively.

When $\tau = \{\theta\}$, the distribution degenerate at $\theta \in \Theta$, we will simply write $d_{Y|X}(\tau, \pi) = d_{Y|X}(\theta, \pi)$, where

$$(2.2) \qquad d_{Y|X}(\theta, \pi) = \iint \log\left\{\frac{p(y|x,\theta)}{p^\pi(y|x)}\right\} p(x,y|\theta)\, d\mu(x)\, d\mu(y)$$

is the expected regret under the loss function $-\log p^\pi(y|x)$ associated with using the predictive density $p^\pi(y|x)$ when $X$ and $Y$ arise from $p(x|\theta)$ and

$p(y|x,\theta)$, respectively. The regret (2.2) is the conditional relative entropy $D(p(Y|X,\theta)\|p^\pi(Y|X))$. The readily derived relationship

$$(2.3) \qquad \int d_{Y|X}(\theta,\pi)\,d\tau(\theta) = d_{Y|X}(\tau,\pi) + \int d_{Y|X}(\theta,\tau)\,d\tau(\theta)$$

implies that (2.2) is a proper scoring rule, as pointed out by Aitchison [1]; that is, the left-hand side of (2.3) attains its minimum value over $\pi \in \Pi$ when $\pi = \tau$. We note that the final integral in (2.3) is the Shannon conditional mutual information $I(Y;\theta|X)$ between $Y$ and $\theta$ conditional on $X$ (under the prior $\tau$). Conditional mutual information has been used by Sun and Berger [21] for deriving reference priors conditional on a parameter to which a subjective prior has been assigned, and by Clarke and Yuan [9] for deriving possibly data-dependent "partial information" reference priors that are conditional on a statistic.

Definition (2.1) of the posterior predictive regret is motivated by standard arguments for adopting the logarithmic score $\log q(Y)$ as an operational utility function when using $q$ as a predictive density for the random quantity $Y$; see, for example, the discussion in Chapter 2 of [6]. The criterion (2.2) was used by Aitchison [1] for the purpose of comparing the predictive performance of estimative and posterior predictive distributions, which was followed up by Komaki [16], who considered the associated asymptotic theory for curved exponential families. Hartigan [14] obtained related higher-order asymptotic expressions which he used to compare estimative predictive distributions based on (bias-corrected) maximum likelihood and Bayes estimators. Akaike [2] discussed the use of (2.2) for the selection of objective priors. A similar approach was also proposed by Geisser in his discussion of Bernardo [5]. Recently, Liang and Barron [19] have derived exact minimax priors under the criterion (2.2) for location and scale families.

The criterion (2.1) extends the domain of definition of (2.2) from degenerate priors $\{\theta\}$ to all proper priors $\tau \in \Omega$. We argue that (2.1) is a suitable Bayesian performance characteristic for assessing the predictive performance of a nonsubjective prior distribution $\pi$ when $\theta$ arises from alternative proper prior distributions $\tau$. There are two ways of thinking about this. First, we might be interested in the predictive performance of a proposed nonsubjective prior distribution under its repeated *use*, as opposed to its performance under repeated *sampling*, as measured by (2.2). From this point of view, we could consider the prior selection problem as an idealized game between the Statistician and Nature, in which each player selects a prior distribution. An alternative viewpoint is to consider (2.1) as measuring the predictive performance of $\pi$ in relation to a subjective prior distribution $\tau$ that is as yet unspecified. Thus, $\tau$ might reflect the prior beliefs, yet to be elicited, of an expert. In this case the prior selection problem could be viewed as a game between the Statistician and an Expert. It is possible, of course, that the Statistician and Expert are the same person, whose prior beliefs have yet to be properly formulated.

Akaike [2] considered priors that give constant posterior predictive regret (2.2), referring to such priors as uniform or "impartial" priors. Such priors will only exist in special cases, however. Achieving constant regret over all possible priors $\tau \in \Omega$ in (2.1) is clearly never possible since, for any fixed $\pi \in \Pi$, the precision of the predictive distribution under $\tau$ will tend to increase as $\tau$ becomes more informative, in which case $d_{Y|X}(\tau, \pi)$ will eventually increase. Alternatively, since $\tau$ is unknown, one might wish to consider the minimaxity of $\pi$ over all $\tau \in \Omega$. However, the maximum regret will tend to occur at degenerate $\tau$. We would therefore be led back to the frequentist risk criterion (2.2), which is not the object of primary interest in the present paper.

For these reasons, we will study the loss function

$$(2.4) \qquad L_{Y|X}(\tau, \pi; \pi^B) = d_{Y|X}(\tau, \pi) - d_{Y|X}(\tau, \pi^B),$$

provided that this exists (see later), which is the posterior predictive regret associated with using the prior $\pi$ compared to using a fixed base prior $\pi^B \in \Pi$. Since we will be investigating default priors for prediction, it is necessary that our procedure for choosing the base measure $\pi^B$ is such that $p^B(y|x)$ does not depend on the particular parameterization of the model that is adopted. We are therefore inevitably led to a choice of base measure that is invariant under arbitrary reparameterization. In the case of a regular parametric family, an obvious candidate for $\pi^B$ is Jeffreys' invariant prior with density proportional to $|I(\theta)|^{1/2}$, where $I(\theta)$ is Fisher's information in the sample $X$. Since we will only be considering regular likelihoods in the rest of this paper, we take $\pi^B = \pi^J$ in the sequel and simply write $L_{Y|X}(\tau, \pi; \pi^J) = L_{Y|X}(\tau, \pi)$.

Assume that the base Jeffreys' prior $\pi^J$ satisfies $d_{Y|X}(\theta, \pi^J) < \infty$ for all $\theta \in \Theta$ and let $p^J(y|x)$ be the conditional density of $Y$ given $X$ under $\pi^J$. Then the (*posterior*) *predictive loss function* defined by

$$(2.5) \qquad \begin{aligned} L_{Y|X}(\theta, \pi) &= d_{Y|X}(\theta, \pi) - d_{Y|X}(\theta, \pi^J) \\ &= \int\!\!\int \log\!\left\{\frac{p^J(y|x)}{p^\pi(y|x)}\right\} p(x, y|\theta) \, d\mu(x) \, d\mu(y) \end{aligned}$$

is well defined, although possibly $+\infty$. Now let $\Omega_{Y|X} \subset \Omega$ be the class of proper priors $\tau$ for which $\int d_{Y|X}(\theta, \pi^J) \, d\tau(\theta) < \infty$. Then for $\pi \in \Pi$ and $\tau \in \Omega_{Y|X}$, we can define the expected predictive loss

$$(2.6) \qquad \begin{aligned} L_{Y|X}(\tau, \pi) &= \int L_{Y|X}(\theta, \pi) \, d\tau(\theta) \\ &= \int d_{Y|X}(\theta, \pi) \, d\tau(\theta) - \int d_{Y|X}(\theta, \pi^J) \, d\tau(\theta) \\ &= d_{Y|X}(\tau, \pi) - d_{Y|X}(\tau, \pi^J), \end{aligned}$$

as in (2.4). Since $\tau \in \Omega_{Y|X}$, the final line is well defined (possibly $+\infty$).

Next we define, for $\tau \in \Omega$,

$$(2.7) \quad \zeta_{Y|X}(\tau) = d_{Y|X}(\tau, \pi^J) = \int\int \log\left\{\frac{p^\tau(y|x)}{p^J(y|x)}\right\} p^\tau(x,y) \, d\mu(x) \, d\mu(y).$$

Since the negative conditional relative entropy $-d_{Y|X}(\tau, \pi^J) = -D(p^\tau(Y|X) \| p^J(Y|X))$ is a natural information-theoretic measure of the uncertainty in the predictive distribution $p^\tau(Y|X)$, we will refer to $\zeta_{Y|X}(\tau)$ as the *predictive information* in $\tau$. Here $p^J(y|x)$ acts as a normalization of the conditional entropy of $p^\tau(y|x)$. From relation (2.3) with $\pi = \pi^J$, we see that $\zeta_{Y|X}(\tau) \leq \int d_{Y|X}(\theta, \pi^J) \, d\tau(\theta)$, from which it follows that $\sup_{\tau \in \Omega} \zeta_{Y|X}(\tau) = \sup_{\theta \in \Theta} \zeta_{Y|X}(\{\theta\})$. That is, the maximum predictive information occurs at (or near) a degenerate prior. Thus, $\zeta_{Y|X}(\tau)$ is a natural entropy-based measure of the information in the predictive distribution $p^\tau(y|x)$. Note that, again from (2.3), $\zeta_{Y|X}(\tau) < \infty$ whenever $\tau \in \Omega_{Y|X}$.

It now follows from (2.3), (2.6) and (2.7) that, for $\pi \in \Pi$ and $\tau \in \Omega_{Y|X}$, we can write

$$(2.8) \qquad d_{Y|X}(\tau, \pi) = L_{Y|X}(\tau, \pi) + \zeta_{Y|X}(\tau).$$

We will explore priors for which $L_{Y|X}(\theta, \pi)$ is approximately constant in $\theta \in \Theta$. Notice that if $L_{Y|X}(\theta, \pi)$ is approximately constant, then, from (2.8), $d_{Y|X}(\tau, \pi)$ is approximately constant over all $\tau$ having the same predictive information $\zeta_{Y|X}(\tau)$. This therefore provides a suitable notion of approximate uniformity of the posterior predictive regret (2.1).

In Sections 4 and 5 we will derive large sample forms, $L(\theta, \pi), L(\tau, \pi), \zeta(\tau)$ and $d(\tau, \pi)$, respectively, of suitably normalized versions of $L_{Y|X}(\theta, \pi)$, $L_{Y|X}(\tau, \pi), \zeta_{Y|X}(\tau)$ and $d_{Y|X}(\tau, \pi)$ and simply refer to $L(\theta, \pi)$ as the *predictive loss function*. Importantly, for smooth priors $\pi$ this asymptotic loss function will not depend on the amount of prediction $Y$ to be carried out. In Section 6 we will investigate uniform and minimax priors under predictive loss. As is often the case in game theory, there is a strong relationship between constant loss, minimax and maximin priors. We give an informal statement of Theorem 6.1. An *equalizer prior* is a prior $\pi$ for which the predictive loss function $L(\theta, \pi)$ is constant over $\theta \in \Theta$. Suppose that $\pi_0$ is an equalizer prior and that there exists a sequence $\tau_k$ of proper priors in the class $\Phi \subset \Omega$, to be defined in Section 4, for which $d(\tau_k, \pi_0) \to 0$ as $k \to \infty$. Then Theorem 6.1 states that $\pi_0$ is minimax with respect to $L(\tau, \pi)$ and $\zeta(\pi_0) = \inf_{\tau \in \Phi} \zeta(\tau)$; that is, $\pi_0$ contains minimum predictive information about $Y$. This latter property is equivalent to $\pi_0$ being maximin, or least favorable, under $L(\tau, \pi)$. Since by construction $L(\tau, \pi^J) = 0$ for all $\tau \in \Phi$, $\pi^J$ is automatically an equalizer prior. However, there may not exist a sequence $\tau_k$ of proper priors with $d(\tau_k, \pi^J) \to 0$, in which case Jeffreys' prior may not be minimax. Some examples will be given in Section 6.

Although the focus of this paper is on the general asymptotic form of the predictive loss, we briefly note the implications of adopting either the posterior predictive regret (2.2) or the predictive loss (2.5) in the special case where the family $p(\cdot|\theta)$ of densities is invariant under a suitable group $\mathcal{G}$ of transformations of the sample space. See, for example, Chapter 6 in [4] for a general discussion of invariant decision problems. Let $\overline{\mathcal{G}}$ be the induced group of transformations on $\Theta$. Then the predictive loss (2.5) is invariant under $\overline{\mathcal{G}}$ and the invariant decisions are invariant priors satisfying $\pi(\bar{g}(\theta)) \propto \pi(\theta)|d\theta/d\bar{g}(\theta)|$ for all $\bar{g} \in \overline{\mathcal{G}}$. If the group $\overline{\mathcal{G}}$ is transitive, then the predictive loss is constant for every invariant prior. Furthermore, if we consider the broader decision problem in which we replace $p^\pi(\cdot|x)$ by the arbitrary decision function $\delta(x) = q_x$, where $q_x(\cdot)$ is to be used as a predictive density for $Y$ when $X = x$, then it can be shown that $p^R(y|x)$, the posterior predictive density under the right Haar measure on $\Theta$, is the best invariant predictive density under the posterior predictive regret (2.2). Since $\pi^J$ is an invariant prior, it further follows that the right Haar measure is the best invariant prior under the predictive loss function (2.5). Since submission of the final version of the present paper, a careful analysis using (2.2) for location and scale families has appeared in [19].

Returning to the definition of the predictive loss function (2.4) relative to an arbitrary base measure $\pi^B$, we see that this is related to the expected predictive loss (2.6) by the equation

$$L_{Y|X}(\tau, \pi; \pi^B) = L_{Y|X}(\tau, \pi) - L_{Y|X}(\tau, \pi^B).$$

Therefore, using $\pi^B$ will give rise to an equivalent predictive loss function if and only if $L_{Y|X}(\theta, \pi^B)$ is constant in $\theta$. In this case we say that $\pi^B$ is *neutral* relative to $\pi^J$.

**3. Relationship to prior predictive regret.** In this section we relate the posterior predictive regret (2.2) and loss function (2.5) to the prior predictive regret and loss function. We will use these relationships in Section 4 to obtain the asymptotic posterior predictive regret $d(\tau, \pi)$ and loss $L(\tau, \pi)$.

For $\pi \in \Pi$, we define the *prior predictive regret* by

$$(3.1) \quad d_X(\theta, \pi) = D(p(X|\theta) \| p^\pi(X)) = \int \log\left\{\frac{p(x|\theta)}{p^\pi(x)}\right\} p(x|\theta) \, d\mu(x),$$

which is the relative entropy $D(p(X|\theta) \| p^\pi(X))$ between $p(x|\theta)$ and the prior predictive density $p^\pi(x)$. Note that $\pi$ may be improper in this definition. In that case, unlike the posterior predictive regret, alternative normalizing constants will give rise to alternative versions of (3.1), differing by constants. The prior predictive regret (3.1) is the focus of work by Bernardo [5], Clarke and Barron [7] and others. Now define $\Pi_X \subset \Pi$ to be the class of priors $\pi$ in $\Pi$ for which $d_X(\theta, \pi) < \infty$ for all $\theta \in \Theta$. If $\pi^J \in \Pi_X$, then for $\pi \in \Pi$ we define the *prior predictive loss* by

$$(3.2) \quad L_X(\theta, \pi) = d_X(\theta, \pi) - d_X(\theta, \pi^J) = \int \log\left\{\frac{p^J(x)}{p^\pi(x)}\right\} p(x|\theta) \, d\mu(x),$$

which is well defined (possibly $+\infty$).

The posterior predictive regret (2.2) and loss (2.5) are simply related to the prior predictive regret (3.1) and loss (3.2). The following result is essentially the chain rule for relative entropy. However, we formally state and prove it since, first, the distribution of $X$ may be improper here and, second, we need to make sure that these relationships are well defined.

LEMMA 3.1. *Suppose that $\pi \in \Pi_{X,Y}$. Then $\pi \in \Pi_X$, $d_{Y|X}(\theta, \pi) < \infty$ for all $\theta \in \Theta$ and*

$$(3.3) \quad d_{Y|X}(\theta, \pi) = d_{X,Y}(\theta, \pi) - d_X(\theta, \pi).$$

*If further $\pi^J \in \Pi_{X,Y}$, then $L_{Y|X}(\theta, \pi) < \infty$ for all $\theta \in \Theta$ and*

$$(3.4) \quad L_{Y|X}(\theta, \pi) = L_{X,Y}(\theta, \pi) - L_X(\theta, \pi).$$

PROOF. Since $\pi \in \Pi$, the marginal densities $p^\pi(X)$ and $p^\pi(X, Y)$ are a.s. ($\theta$) finite for all $\theta \in \Theta$. Therefore,

$$p^\pi(x,y) = \int p(x,y|\phi) \, d\pi(\phi) = p^\pi(x) \int p(y|x,\phi) \, dp^\pi(\phi|x) = p^\pi(x) p^\pi(y|x),$$

since, by definition, $p(x|\phi) \, d\pi(\phi) = p^\pi(x) \, dp^\pi(\phi|x)$. It now follows straightforwardly from the definitions (2.2) and (3.1) that

$$(3.5) \quad d_{X,Y}(\theta, \pi) = d_{Y|X}(\theta, \pi) + d_X(\theta, \pi).$$

Since $\pi \in \Pi_{X,Y}$, it follows from (3.5) that both $d_{Y|X}(\theta, \pi) < \infty$ and $\pi \in \Pi_X$ and, hence, relation (3.3) holds. Since $\pi \in \Pi_X$ and $\pi^J \in \Pi_X$, it follows from (3.2) that $L_{Y|X}(\theta, \pi)$ is finite for all $\theta$. Finally, since $\pi^J \in \Pi$, we have $p^J(x, y) = p^J(x) p^J(y|x)$ and relation (3.4) follows straightforwardly from the definitions (2.5) and (3.2).

Finally, let $\Omega_X \subset \Omega$ be the class of priors $\tau$ in $\Omega$ satisfying $\int d_X(\theta, \pi^J) \, d\tau(\theta) < \infty$. It follows from equation (3.3) of Lemma 3.1 that $\pi^J \in \Pi_{X,Y}$ and $\tau \in \Omega_{X,Y}$ imply that $\int d_{Y|X}(\theta, \pi^J) \, d\tau(\theta) < \infty$, $\tau \in \Omega_X$ and

$$\int d_{Y|X}(\theta, \pi^J) \, d\tau(\theta) = \int d_{X,Y}(\theta, \pi^J) \, d\tau(\theta) - \int d_X(\theta, \pi^J) \, d\tau(\theta).$$

Therefore, if $\pi \in \Pi_{X,Y}$ and $\tau \in \Omega_{X,Y}$, then the expected posterior loss $L_{Y|X}(\tau, \pi)$ at (2.6) is well defined. $\square$

**4. Asymptotic behavior of the predictive loss.** Throughout the remainder of this article we specialize to the case $X = (X_1, \ldots, X_n)$ and $Y = (X_{n+1}, \ldots, X_{n+m})$, where the $X_i$ are independent observations from a density $f(x|\theta)$ with respect to a measure $\mu$. In the present section we investigate the asymptotic behavior as $n \to \infty$ of the predictive loss function (2.5). In particular, we will show that, under suitable regularity conditions, the asymptotic form of (2.5) (after suitable normalization) is the same regardless of the amount $m$ of prediction to be performed. This leads to a general definition for broad classes of priors $\pi$ and $\tau$ of the (asymptotic) predictive loss $L(\tau, \pi)$, information $\zeta(\tau)$ and regret $d(\tau, \pi)$.

For an asymptotic analysis of the posterior predictive regret (2.2) and loss function (2.5), from (3.2), (3.3) and (3.4), we see that it suffices to study the asymptotic behavior of the prior predictive regret $d_X(\theta, \pi)$. Suppose that $\pi \in \Pi$ has a density with respect to Lebesgue measure. For notational convenience, in what follows we will use the same symbol $\pi$ to denote this density. Let $l(\theta) = n^{-1} \log p(X|\theta) = n^{-1} \sum_{i=1}^n \log f(X_i|\theta)$ be the normalized loglikelihood function and let $i(\theta) = E^\theta\{-l''(\theta)\} = n^{-1} I(\theta)$ be Fisher's information per observation. A standard result for the prior predictive regret (3.1) when $\pi$ is a density (see, e.g., [7]) is that, under suitable regularity conditions,

$$(4.1) \qquad d_X(\theta, \pi) = \frac{p}{2} \log\left(\frac{n}{2\pi e}\right) + \log\left\{\frac{|i(\theta)|^{1/2}}{\pi(\theta)}\right\} + o(1)$$

as $n \to \infty$. [Here the $\pi$ appearing in the first term on the right-hand side of (4.1) is the usual transcendental number and should not be confused with the prior $\pi(\cdot)$.] Taking Jeffreys' prior to be $\pi^J(\theta) = |i(\theta)|^{1/2}$, it follows from (3.2) and (4.1) that the prior predictive loss satisfies

$$L_X(\theta, \pi) = \log\left\{\frac{|i(\theta)|^{1/2}}{\pi(\theta)}\right\} + o(1).$$

It now follows from (3.4) that, for any sequence $m = m_n \geq 1$, $L_{Y|X}(\theta, \pi) = o(1)$; that is, to first order the posterior predictive loss is identically zero for every smooth prior $\pi$. It is therefore necessary to develop further the asymptotic expansion in (4.1). Let $\hat{\theta}$ denote the maximum likelihood estimator based on the data $X$ and assume that the observed information matrix $J = -nl''(\hat{\theta})$ is positive definite over the set $S$ for which $P^\theta(S) = 1 + o(n^{-1})$, uniformly in compact subsets of $\Theta$.

Let $\Pi_\infty$ be the class of priors $\pi \in \Pi$ for which $\pi \in \Pi_X$ for all $n$ and let $C \subset \Pi_\infty$ be the class of priors in $\Pi_\infty$ that possess densities having continuous second-order derivatives throughout $\Theta$. Then, under suitable additional regularity conditions on $f$ and $\pi \in C$ to be discussed in Section 5, the marginal density of $X$ is

$$p^\pi(x) = (2\pi s_B^2)^{p/2} |J|^{-1/2} p(x|\hat{\theta}) \pi(\hat{\theta}) \{1 + o(n^{-1})\},$$

where $s_B^2 = (1 + b_B)^2$ is a *Bayesian Bartlett correction*, with $b_B = O(n^{-1})$; see, for example, [22]. Therefore, we can write

$$\log\left\{\frac{p(x|\theta)}{p^\pi(x)}\right\} = \frac{p}{2}\log\left(\frac{n}{2\pi e}\right) + \log\left\{\frac{|i(\theta)|^{1/2}}{\pi(\theta)}\right\} - pb_B - \left[n\{l(\hat{\theta}) - l(\theta)\} - \frac{p}{2}\right]$$
$$- \log\left\{\frac{\pi(\hat{\theta})}{\pi(\theta)}\right\} + \frac{1}{2}\log\left\{\frac{|J|}{|I(\theta)|}\right\} + o\left(\frac{1}{n}\right).$$

Since $E^\theta[n\{l(\hat{\theta}) - l(\theta)\}] = ps_F^2(\theta)/2 + o(n^{-1})$, where $s_F^2(\theta) = \{1 + b_F(\theta)\}^2$ is a *frequentist Bartlett correction*, with $b_F(\theta) = O(n^{-1})$, it follows from (3.1) that

$$(4.2) \qquad d_X(\theta, \pi) = \frac{p}{2}\log\left(\frac{n}{2\pi e}\right) + \log\left\{\frac{|i(\theta)|^{1/2}}{\pi(\theta)}\right\} - h_n(\theta, \pi),$$

where

$$(4.3) \qquad \begin{aligned} h_n(\theta, \pi) &= p\{E^\theta(b_B) + b_F(\theta)\} + E^\theta\left[\log\left\{\frac{\pi(\hat{\theta})}{\pi(\theta)}\right\}\right] \\ &\quad - \frac{1}{2}E^\theta\left[\log\left\{\frac{|J|}{|I(\theta)|}\right\}\right] + o\left(\frac{1}{n}\right). \end{aligned}$$

Under suitable regularity conditions, the leading term in (4.3) turns out to be $O(n^{-1})$, since both the Bayesian and frequentist Bartlett corrections are $O(n^{-1})$, as are all the expectations on the right-hand side of (4.3). We will therefore suppose that $h_n$ is of the form

$$(4.4) \qquad h_n(\theta, \pi) = \left\{\frac{D(\theta, \pi)}{2n}\right\} + r_n(\theta, \pi),$$

where $D(\theta, \pi)$ is continuous in $\theta$ and the remainder term $r_n(\theta, \pi)$ satisfies one of the following three successively stronger conditions:

R1. $r_n(\theta, \pi) = o(n^{-1})$ uniformly in compacts of $\Theta$;
R2. $r_n(\theta, \pi) = O(n^{-2})$ uniformly in compacts of $\Theta$;
R3. $r_n(\theta, \pi) = E(\theta, \pi)n^{-2} + o(n^{-2})$ uniformly in compacts of $\Theta$, where $E(\theta, \pi)$ is continuous in $\theta$.

The above three forms of remainder require successively stronger assumptions about both the likelihood $p(\cdot|\theta)$ and the prior $\pi(\theta)$. Suitable sets of regularity conditions for the validity of (4.4) will be discussed in Section 5. In particular, $\pi \in C$ is a sufficient condition on the prior for the weakest form R1 of remainder. The form of $D(\theta, \pi)$ for $\pi \in C$ will be derived in Section 5.

Throughout the remainder of the paper we assume that $\pi^J \in C$ and define, for all $\pi \in C$,

$$(4.5) \qquad L(\theta, \pi) = D(\theta, \pi) - D(\theta, \pi^J).$$

We note that $L(\theta, \pi)$ is well defined when $\pi$ is improper since the arbitrary normalizing constant in $\pi$ does not appear in $D(\theta, \pi)$. We will study the asymptotic behavior of the posterior predictive loss (2.5) as $n \to \infty$ for an arbitrary number $m_n \geq 1$ of predictions $Y_i$. Let $c_n = 2n(n + m_n)/m_n$. The next theorem gives conditions under which

$$(4.6) \qquad c_n L_{Y|X}(\theta, \pi) \to L(\theta, \pi)$$

uniformly in compacts of $\Theta$ under each of the forms R1–R3 of remainder.

THEOREM 4.1.

(a) *Suppose that* R1 *holds. Then* (4.6) *holds whenever* $\liminf_{n \to \infty} m_n/n > 0$.

(b) *Suppose that* R2 *holds. Then* (4.6) *holds whenever* $m_n \to \infty$.

(c) *Suppose that* R3 *holds. Then* (4.6) *holds for every sequence* $(m_n)$ *of positive integers.*

PROOF. First note that (3.2), (4.2), (4.4) and (4.5) give, on taking $\pi^J(\theta) = |i(\theta)|^{1/2}$,

$$(4.7) \qquad L_X(\theta, \pi) = \log\left\{\frac{|i(\theta)|^{1/2}}{\pi(\theta)}\right\} - \left\{\frac{L(\theta, \pi)}{2n}\right\} - \bar{r}_n(\theta, \pi),$$

where $\bar{r}_n(\theta, \pi) = r_n(\theta, \pi) - r_n(\theta, \pi^J)$. Also note that, since $\pi \in \Pi_\infty$, Lemma 3.1 applies for all $n$.

(a) From (3.4), (4.7) and R1, we have $L_{Y|X}(\theta, \pi) = c_n^{-1} L(\theta, \pi) + o(n^{-1})$ and (4.6) follows since $n^{-1} c_n = 2(m_n^{-1} n + 1)$ and $\limsup_{n \to \infty} m_n^{-1} n < \infty$.

(b) From (3.4), (4.7) and R2, we have $L_{Y|X}(\theta, \pi) = c_n^{-1} L(\theta, \pi) + O(n^{-2})$ and (4.6) follows since $n^{-2} c_n = 2(m_n^{-1} + n^{-1}) \to 0$.

(c) From (3.4), (4.7) and R3, we have $L_{Y|X}(\theta, \pi) = c_n^{-1}\{L(\theta, \pi) + d_n^{-1} \overline{E}(\theta, \pi)\} + o(n^{-2})$, where $d_n = \{2(2n + m_n)\}^{-1} n(n + m_n)$ and $\overline{E}(\theta, \pi) = E(\theta, \pi) - E(\theta, \pi^J)$. (4.6) follows since $d_n^{-1} = O(n^{-1})$ and $n^{-2} c_n = 2(m_n^{-1} + n^{-1})$ is bounded. □

Theorem 4.1 tells us that, although the predictive loss function (2.5) covers an infinite variety of possibilities for the amount of data to be observed and predictions to be made, it is approximately equivalent to the single loss function (4.5), provided that a sufficient amount of data $X$ is to be observed. Although this is not surprising given the form of (4.7) and the relation (3.4), it considerably simplifies the task of assessing the predictive risk arising from using alternative priors. We will refer to $L(\theta, \pi)$ as the (*asymptotic*) *predictive loss function*. A special case of interest arises when $m_n = n$, which corresponds to prediction of a replicate data set of the same size as that to be observed. Note that in this case (4.6) holds under the weakest condition R1. More generally, Laud and Ibrahim [18] refer to the posterior predictive density of $Y$ in this case as the "predictive density of a replicate experiment," which they study in relation to model choice.

Now let $\Omega_\infty$ be the class of priors $\tau \in \Omega$ for which $\tau \in \Omega_X$ for all $n$. Although the expected predictive loss $L_{Y|X}(\tau, \pi)$ is well defined (possibly $+\infty$) when $\pi \in \Pi_\infty$ and $\tau \in \Omega_\infty$, in general, the expected asymptotic predictive loss $\int L(\theta, \pi) \, d\tau(\theta)$ may not exist, and when it does, additional conditions will be needed for it to be the limit of the expected loss $c_n L_{Y|X}(\tau, \pi)$. In order to retain generality, we will extend the domain of definition of the asymptotic predictive loss (4.5) so that it is defined for all $\pi \in \Pi_\infty$ and $\tau \in \Omega_\infty$. Thus, for $\pi \in \Pi_\infty, \tau \in \Omega_\infty$ and a given sequence $(m_n)$ of positive integers, we define the (*asymptotic*) *predictive loss* to be

$$(4.8) \qquad L(\tau, \pi) = \limsup_{n \to \infty} c_n L_{Y|X}(\tau, \pi),$$

which always exists (possibly $+\infty$). Thus, $L(\tau, \pi)$ represents the asymptotically worst-case predictive loss when the prior $\pi$ is used in relation to the alternative proper prior $\tau$. Since the degenerate prior $\tau = \{\theta\}$ is in $\Omega_\infty$, (4.8) also provides a definition of $L(\theta, \pi)$ for all $\pi \in \Pi_\infty, \theta \in \Theta$, which agrees with (4.5) whenever $\pi \in C \subset \Pi_\infty$ and one of the conditions R1–R3 holds.

Now define the (*asymptotic*) *predictive information* contained in $\tau \in \Omega_\infty \cap \Pi_\infty$ to be

$$(4.9) \qquad \zeta(\tau) = -L(\tau, \tau) = \liminf_{n \to \infty} c_n \zeta_{Y|X}(\tau)$$

and let $\Phi \subset \Omega_\infty \cap \Pi_\infty$ be the class of $\tau$ for which $\zeta(\tau) < \infty$. Finally, for $\pi \in \Pi_\infty$ and $\tau \in \Phi$, define

$$(4.10) \qquad d(\tau, \pi) = L(\tau, \pi) + \zeta(\tau),$$

which is the asymptotic form of equation (2.8). The next lemma implies that the predictive loss function (4.8) is a $\Phi$-proper scoring rule and that $d(\tau, \pi)$ is the regret associated with $L(\tau, \pi)$.

LEMMA 4.1. *For all* $\tau \in \Phi$,

$$\inf_{\pi \in \Pi_\infty} L(\tau, \pi) = L(\tau, \tau) = -\zeta(\tau).$$

PROOF. Let $\tau \in \Phi$. By construction, $d(\tau, \tau) = 0$, so we only need to show that $d(\tau, \pi) \geq 0$ for all $\pi \in \Pi_\infty$. Since $\pi \in \Pi_\infty$ and $\tau \in \Omega_\infty \cap \Pi_\infty$, we have $\pi \in \Pi_{X,Y}$ and $\tau \in \Omega_{X,Y} \cap \Pi_{X,Y}$ for all $n$ and, hence, the quantities $L_{Y|X}(\tau, \pi)$ and $L_{Y|X}(\tau, \tau)$ are both well defined. But $L_{Y|X}(\tau, \tau) \leq L_{Y|X}(\tau, \pi)$ and multiplying both sides of this inequality by $c_n$ and taking the $\limsup_{n \to \infty}$ on both sides of the resulting inequality gives $L(\tau,\tau) \leq L(\tau,\pi)$. The result follows from the definition of $d(\tau,\pi)$. $\square$

When $\pi \in C$, $L(\theta,\pi)$ is independent of the sequence $m_n$. In general, however, both $L(\tau,\pi)$ and $\zeta(\tau)$ may depend on the particular sequence $(m_n)$, although we have suppressed this dependence in the notation. Nevertheless, the minimax results of Section 6 will be independent of $(m_n)$.

**5. Derivation of the asymptotic predictive loss function.** In this section we obtain the form of the function $D(\theta,\pi)$ arising in the $O(n^{-1})$ term in the asymptotic expansion of the prior predictive regret $d_X(\theta,\pi)$. This then leads to an expression for the asymptotic predictive loss function $L(\theta,\pi)$ for all $\pi \in C$ via relation (4.5). The computations involved in the determination of $D(\theta,\pi)$, which are similar in nature to computations in [14], are technically quite demanding. Finally, we deduce expressions for the asymptotic posterior predictive regret (4.10) and predictive information (4.9) under certain conditions.

Theorem 5.1 below is the central result of this section. Write $D_j = \partial/\partial \theta^j$, $j = 1, \ldots, p$. Let $\rho = \rho(\theta) = \log \pi(\theta)$ and write $\rho_r = D_r \rho$. We use the summation convention throughout.

THEOREM 5.1. *Assume that one of the conditions* R1–R3 *holds. Then*

$$D(\theta,\pi) = A(\theta,\pi) + M(\theta), \tag{5.1}$$

*where*

$$A(\theta,\pi) = i^{rs}\rho_r\rho_s + 2D_s(i^{rs}\rho_r) \tag{5.2}$$

*and $M(\theta)$ is independent of $\pi$.*

We will prove Theorem 5.1 via four lemmas, each of which evaluates the leading term in one of the terms on the right-hand side of equation (4.3). We discuss suitable sets of regularity conditions following the proof.

For $1 \leq j,k,r,\ldots \leq p$, define $D_{jkr\cdots} = \frac{\partial}{\partial \theta^j}\frac{\partial}{\partial \theta^k}\frac{\partial}{\partial \theta^r}\cdots$, $a_{jkr\cdots} = \{D_{jkr\cdots}l(\theta)\}_{\theta=\hat{\theta}}$, $c_{jr} = -a_{jr}$, $C = (c_{jr})$, $C^{-1} = (c^{jr})$, $\rho_{jk} = D_{jk}\rho$, $\hat{\rho}_{jk\cdots} = \rho_{jk\cdots}(\hat{\theta})$ and

$$k_{jkl\cdots,rst\cdots} = k_{jkl\cdots,rst\cdots}(\theta) = E^\theta\{D_{jkl\cdots}\log f(X_i;\theta)D_{rst\cdots}\log f(X_i;\theta)\}.$$

Also define

$$k_1^* = i^{jr}(\rho_{jr} + \rho_j\rho_r), \qquad k_2^* = 3k_{jrsu}i^{jr}i^{su},$$
$$k_3^* = 3k_{ijr}\rho_s i^{ij}i^{rs}, \qquad k_4^* = 15k_{jrs}k_{uvw}i^{jr}i^{su}i^{vw}$$

and

$$Q_1 = D_{rs}i^{rs}, \qquad Q_2 = k_2^*, \qquad Q_3 = 3D_s(k_{ijr}i^{ij}i^{rs}), \qquad Q_4 = k_4^*.$$

LEMMA 5.1.
$$nE^\theta(b_B) \to \frac{1}{2p}\left(k_1^* + \frac{1}{12}k_2^* + \frac{1}{3}k_3^* + \frac{1}{36}k_4^*\right).$$

PROOF. Comparing with the Bayesian Bartlett correction factor as given in equation (2.6) of [13], we obtain

(5.3) $$b_B = \frac{1}{2pn}\left(H_1 + \frac{1}{12}H_2 + \frac{1}{3}H_3 + \frac{1}{36}H_4\right) + o(n^{-1}),$$

where
$$H_1 = c^{jr}(\hat{\rho}_{jr} + \hat{\rho}_j\hat{\rho}_r), \qquad H_2 = 3a_{jrsu}c^{jr}c^{su},$$
$$H_3 = 3a_{ijr}\hat{\rho}_s c^{ij}c^{rs}, \qquad H_4 = 15a_{jrs}a_{uvw}c^{jr}c^{su}c^{vw}.$$

Noting that $E^\theta(H_a) = k_a^* + o(1), a = 1,\ldots,4$, the lemma follows from (5.3). $\square$

LEMMA 5.2.
$$nb_F(\theta) \to \frac{1}{2p}\left(Q_1 + \frac{1}{12}Q_2 - \frac{1}{3}Q_3 + \frac{1}{36}Q_4\right).$$

PROOF. Comparing with the frequentist Bartlett correction factor as given in equation (2.10) of [13], we obtain

$$b_F(\theta) = \frac{1}{2pn}\left(Q_1 + \frac{1}{12}Q_2 - \frac{1}{3}Q_3 + \frac{1}{36}Q_4\right) + o(n^{-1}),$$

from which the result follows. $\square$

LEMMA 5.3.
$$nE^\theta\left[\log\left\{\frac{\pi(\hat{\theta})}{\pi(\theta)}\right\}\right] \to \rho_r b^r + \frac{1}{2}i^{jr}\rho_{jr},$$

where $b^r = i^{jr}i^{kt}k_{jk,t} + \frac{1}{2}i^{jr}i^{kt}k_{jkt}$.

PROOF. From [20], page 209, we see that

(5.4) $$E^\theta(\hat{\theta}^r) = \theta^r + n^{-1}b^r + o(n^{-1}),$$

(5.5) $$\text{Cov}^\theta(\hat{\theta}^r, \hat{\theta}^s) = n^{-1}i^{rs} + o(n^{-1}).$$

By applying Bartlett's identity,
$$k_{jkt} + k_{j,kt} + k_{k,jt} + k_{t,jk} + k_{j,k,t} = 0$$

(cf. equation (7.2) of [20]), it can be seen that our expression for $b^r$ agrees with that of McCullagh. From (5.4), (5.5) and the Taylor expansion of $\rho(\hat{\theta})$ around $\theta$, we obtain

$$E^\theta\{\rho(\hat{\theta})\} = \rho(\theta) + n^{-1}b^r\rho_r + \tfrac{1}{2}n^{-1}\rho_{rs}i^{rs} + o(n^{-1}),$$

from which the lemma follows. □

LEMMA 5.4.
$$nE^\theta\left[\log\left\{\frac{|J|}{|I(\theta)|}\right\}\right]$$
$$\to -i^{jr}\left(k_{jrs}b^s + i^{sk}k_{jrs,k} + \frac{1}{2}k_{jrst}i^{st}\right)$$
$$-\frac{1}{2}i^{jl}i^{vi}\{(k_{ji,lv} - i_{ji}i_{lv}) + k_{jis}i^{ts}k_{lv,t} + k_{lvw}i^{tw}k_{ji,t} + k_{jit}k_{lvw}i^{tw}\}.$$

PROOF. By the Taylor expansion of $a_{jr} = l_{jr}(\hat{\theta})$ around $\theta$, we get

(5.6) $$a_{jr} = k_{jr}(\theta) + e_{jr} + o(n^{-1}),$$

where

(5.7) $$\begin{aligned}e_{jr} &= l_{jr} - k_{jr} + k_{jrs}(\hat{\theta}^s - \theta^s) \\ &\quad + (l_{jrs} - k_{jrs})(\hat{\theta}^s - \theta^s) + \tfrac{1}{2}k_{jrst}(\hat{\theta}^s - \theta^s)(\hat{\theta}^t - \theta^t).\end{aligned}$$

From (5.6) and (5.7), we obtain
$$C = i(\theta) - E_* + o(n^{-1}),$$

where $E_* = (e_{jr})$. Noting that $J = nC$, $I(\theta) = ni(\theta)$, $i(\theta)$ positive definite and $E_*$ is a matrix with elements of order $O(n^{-1/2})$, from the above expression for $C$ and standard results on the eigenvalues and determinant of a matrix, it follows by the Taylor expansion that

(5.8) $$\log\left\{\frac{|J|}{|I(\theta)|}\right\} = -\operatorname{tr}\{i^{-1}(\theta)E_*\} - \frac{1}{2}\operatorname{tr}\{i^{-1}(\theta)E_*i^{-1}(\theta)E_*\} + o(n^{-1/2}).$$

Using an expansion for $\hat{\theta}^s - \theta^s$ as in [20], Chapter 7, we obtain

(5.9) $$\hat{\theta}^s - \theta^s = i^{js}\{l_j + i^{uk}l_u(l_{jk} - k_{jk}) + \tfrac{1}{2}k_{jkt}i^{uk}i^{wt}l_ul_w\} + o(n^{-1/2}).$$

Substituting (5.9) into (5.7) and using (5.4) and (5.5), it follows that

(5.10) $$E^\theta(e_{jr}) = n^{-1}(k_{jrs}b^s + k_{jrs,k}i^{sk} + \tfrac{1}{2}k_{jrst}i^{st}) + o(n^{-1})$$

and

(5.11) $$\begin{aligned}E^\theta(e_{jr}e_{ku}) = n^{-1}\{&(k_{jr,ku} - i_{jr}i_{ku}) \\ &+ (k_{jrt}k_{ku,w} + k_{kuw}k_{jr,t} + k_{jrt}k_{kuw})i^{tw}\} + o(n^{-1}).\end{aligned}$$

While all four terms on the right-hand side of (5.7) are required in evaluating (5.10), only the first two terms on the right-hand side of (5.7) are required in evaluating (5.11). The lemma follows on taking expectations on both sides of (5.8) and using (5.10) and (5.11) on the right-hand side. □

PROOF OF THEOREM 5.1. First, putting Lemmas 5.1 and 5.2 together gives

$$np\{E^\theta(b_B) + b_F(\theta)\} \to \tfrac{1}{2}\{(Q_1 + k_1^*) - \tfrac{1}{3}(Q_3 - k_3^*) + \tfrac{1}{6}(Q_2 + \tfrac{1}{3}Q_4)\}.$$

Along with Lemmas 5.3 and 5.4, this gives equation (5.1) with

$$A(\theta, \pi) = i^{rs}(\rho_r \rho_s + 2\rho_{rs}) + 2(k_{jku} + k_{jk,u})i^{ku}i^{jr}\rho_r.$$

Now note that $D_r i_{kj} = -D_r E^\theta(l_{kj}) = -(k_{kjr} + k_{kj,r})$ so that

$$A(\theta, \pi) = i^{rs}(\rho_r \rho_s + 2\rho_{rs}) - 2D_u(i_{jk})i^{ku}i^{jr}\rho_r.$$

Finally, $D_u(i_{jk})i^{ku}i^{jr} = -D_u(i^{ku})i_{jk}i^{jr} = -D_u(i^{ru})$ and so

$$A(\theta, \pi) = i^{rs}(\rho_r \rho_s + 2\rho_{rs}) + 2D_s(i^{rs})\rho_r = i^{rs}\rho_r \rho_s + 2D_s(i^{rs}\rho_r),$$

as required. □

We briefly discuss suitable regularity conditions on the likelihood and prior for the validity of the three forms of remainder R1–R3, although we will not dwell on alternative sets of sufficient conditions in the present paper. There are broadly two sets of conditions required, those for the validity of the Laplace approximation of $p^\pi(x)$ and those for the validity of the approximation of each of the terms in (4.3). Consider first the form of remainder R2, ignoring for the moment the uniformity requirement. A suitable set of conditions for this form of remainder is given in Section 3 of [15], which constitutes the definition of a "Laplace-regular" family. Broadly, one requires $l(\theta)$ to be six-times continuously differentiable and $\pi(\theta)$ to be four-times continuously differentiable, plus additional conditions controlling the error term and nonlocal behavior of the integrand. Since additionally we require uniformity in compact subsets of $\Theta$ in R2, we need to replace the neighborhood $B_\varepsilon(\theta_0)$ in these conditions by an arbitrary compact subset of $\Theta$. In addition to these conditions, for the approximation of the terms in (4.3) we require the expectations of the mixed fourth-order partial derivatives of $\log f(X; \theta)$ to be continuous and also conditions guaranteeing the expansions for the expectation of $\hat\theta$ needed in the proofs of Lemmas 5.3 and 5.4, as given in [20], Chapter 7. From an examination of the relevant proofs, it is seen that a slight strengthening of the above conditions will be required for the stronger form R3 of remainder. For example, $l(\theta)$ and $\pi(\theta)$ seven-times and five-times continuously differentiable, respectively, will give rise to a higher-order version of Laplace-regularity. Finally, the weaker form of remainder R1 would apply when $l(\theta)$ and $\pi(\theta)$ are only four-times and twice continuously differentiable, respectively, again with additional regularity conditions controlling, for example, the nonlocal behavior of the integrand in the Laplace approximation and giving uniformity of all the $o(n^{-1})$ remainder terms.

Returning to the predictive loss function, it follows from Theorem 5.1 that, for $\pi \in C$, the asymptotic predictive loss function (4.5) is given by

$$(5.12) \qquad L(\theta, \pi) = A(\theta, \pi) - A(\theta, \pi^J),$$

where $A(\theta, \pi^J) = i^{rs}\nu_r\nu_s + 2D_s(i^{rs}\nu_r)$ and $\nu = \log \pi^J = \frac{1}{2}\log|i|$. It is interesting to note that (5.12) is of the same form as the right-hand side of the first expression in Theorem 4 of [14], which relates to the comparison of estimative predictive distributions based on Bayes estimators. In the case of a single prediction ($m=1$), the connection can be understood from Theorem 7 of [14], which establishes that, to the asymptotic order considered here, the Kullback–Leibler difference between the posterior and the associated estimative predictive distributions is independent of the prior. The derivation of Theorem 5.1 given here is more direct, as it does not involve Bayes estimators. Moreover, our result applies for an arbitrary amount of prediction.

Note that $L(\theta, \pi)$ only depends on the sampling model through Fisher's information. The quantity $M(\theta)$, however, involves components of skewness and curvature of the model. We do not consider $M(\theta)$ further in this paper, although its form, which may be deduced from the results of Lemmas 5.1–5.4, may be of independent interest. It may be verified directly that $L(\theta, \pi)$ is invariant under parameter transformation, as expected in view of (4.6) and the invariance of $L_{Y|X}(\theta, \pi)$. Furthermore, since all the terms in (4.2) are invariant, it follows that $\overline{M}(\theta) \equiv M(\theta) + A(\theta, \pi^J)$ must also be an invariant quantity. In the case $p=1$, we obtain the relatively simple expression

$$(5.13) \qquad \overline{M}(\theta) = \tfrac{1}{12}\alpha_{111}^2 + \tfrac{1}{2}\gamma^2,$$

where $\alpha_{111}$ is the skewness and $\gamma^2 = \alpha_{22} - \alpha_{12}^2 - 1$ is Efron's curvature, with

$$\alpha_{jk\cdots}(\theta) = \{i(\theta)\}^{-(j+k+\cdots)/2} E^\theta\{l^j(\theta)l^k(\theta)\ldots\},$$

where $l^j$ is the $j$th derivative of $l$.

EXAMPLE 5.1. Normal model with unknown mean. As a simple first example, suppose that $X_i \sim N(\theta, 1)$. Here $i(\theta) = 1$ and $\alpha_{111}(\theta) = \gamma^2(\theta) = 0$ so that $L(\theta, \pi) = (\rho')^2 + 2\rho''$ and $\overline{M}(\theta) = 0$ from (5.13). By construction, $L(\theta, \pi^J) = 0$, but note that the improper priors $\pi^c \propto \exp\{c(\theta - \theta_0)\}, c \in \mathcal{R}$, also deliver constant loss, with $L(\theta, \pi^c) = c^2 > 0$. We will see in Section 6 that Jeffreys' prior is minimax in this example. Since here $\overline{M}(\theta) = 0$ and $\pi^J(\theta) \propto 1$, this result also follows from the exact analysis of the criterion (2.1) in [19].

Now let $\overline{\Omega}$ be the class of priors having compact support in $\Theta$ and let $\Gamma = \overline{\Omega} \cap C$. It follows from (4.6) that if $\pi \in C$ and $\tau \in \overline{\Omega}$, then $L(\tau, \pi)$ is equal to the expected predictive loss $\int L(\theta, \pi) \tau(\theta) \, d\theta$. Since $\tau \in C$, we also have $\zeta(\tau) = -\int L(\theta, \tau) \tau(\theta) \, d\theta$, which is finite since $L(\theta, \tau)$ is continuous and, hence, bounded on compact subsets of $\Theta$. The next result gives expressions for the predictive regret $d(\tau, \pi)$ and predictive information $\zeta(\tau)$ when $\pi \in C$ and $\tau \in \Gamma$. The expression for $\zeta(\tau)$ here is similar to that given in Theorem 5 of [14] for the Bayes risk of bias-adjusted estimators.

LEMMA 5.5. *Suppose $\pi \in C$ and $\tau \in \Gamma$. Then*

$$(5.14) \qquad d(\tau, \pi) = \int i^{rs} (\rho_r - \mu_r)(\rho_s - \mu_s) \tau \, d\theta$$

*and*

$$(5.15) \qquad \zeta(\tau) = \int i^{rs} (\mu_r - \nu_r)(\mu_s - \nu_s) \tau \, d\theta,$$

*where $\mu = \log \tau$.*

PROOF. From (5.2), integration by parts gives

$$(5.16) \quad \int A(\theta, \pi) \tau(\theta) \, d\theta = \int i^{rs} \rho_r \rho_s \tau \, d\theta - 2 \int i^{rs} \rho_r \mu_s \tau \, d\theta + 2\beta(\tau, \pi),$$

where

$$\beta(\tau, \pi) = \sum_{s=1}^{p} \int [i^{rs} \rho_r \tau]_{\underline{\theta}^s(\theta^{(-s)})}^{\bar{\theta}^s(\theta^{(-s)})} \, d\theta^{(-s)}$$

and $\underline{\theta}^s(\theta^{(-s)})$ and $\bar{\theta}^s(\theta^{(-s)})$ are the finite lower and upper limits of integration for $\theta^s$ for fixed $\theta^{(-s)}$, the vector of components of $\theta$ omitting $\theta^s$. But $\beta(\tau, \pi) = 0$, since both $\pi$ and $\tau$ are in $C$. Therefore,

$$(5.17) \qquad \int A(\theta, \pi) \tau(\theta) \, d\theta = \int i^{rs} \rho_r (\rho_s - 2\mu_s) \tau \, d\theta.$$

Evaluating (5.17) at $\pi = \tau \in C$ gives

$$(5.18) \qquad \int A(\theta, \tau) \tau(\theta) \, d\theta = -\int i^{rs} \mu_r \mu_s \tau \, d\theta.$$

It now follows from (5.17) and (5.18) that

$$d(\tau, \pi) = L(\tau, \pi) - L(\tau, \tau) = \int \{A(\theta, \pi) - A(\theta, \tau)\} \tau(\theta) \, d\theta$$

$$= \int i^{rs} \{\rho_r(\rho_s - 2\mu_s) + \mu_r \mu_s\} \tau \, d\theta,$$

which gives (5.14). Since $\zeta(\tau) = d(\tau, \pi^J)$, (5.15) follows on evaluating the above expression at $\pi = \pi^J$. $\square$

The expression (5.15) for the predictive information $\zeta(\tau)$ is seen to be invariant under reparameterization, as expected. It might appear at first sight that $\zeta(\tau)$ will attain the value zero at $\tau = \pi^J$, but this is not necessarily the case since $\pi^J$ may be improper and there may be no sequence of priors in $\Gamma$ converging to $\pi^J$ in the right way: see the next section. Finally, note that the form of $d(\tau, \pi)$ in Lemma 5.5 implies that $L(\theta, \pi)$ is a $\Gamma$-strictly proper scoring rule since $d(\tau, \pi)$ attains its minimum value of zero uniquely at $\pi = \tau \in \Gamma$.

**6. Impartial, minimax and maximin priors.** As expected, for a given prior density $\pi \in \Pi_\infty$, from (4.10) the posterior predictive regret will be large when the predictive information (4.9) in $\tau$ is large. Therefore it is not possible to achieve constant regret over all possible $\tau \in \Phi$, nor minimaxity since the regret is unbounded. Instead, as discussed in Section 2, we consider the predictive regret associated with using $\pi$ compared to using Jeffreys' prior and study the behavior of the predictive loss function

$$(6.1) \qquad L(\tau, \pi) = d(\tau, \pi) - d(\tau, \pi^J),$$

which is the asymptotic form of the normalized version of equation (2.4).

Adopting standard game-theoretic terminology, the prior $\pi \in \Pi_\infty$ is an *equalizer* prior if the predictive loss $L(\theta, \pi)$ is constant over $\theta \in \Theta$. This is equivalent to the predictive loss (6.1) being constant over all $\tau \in \Gamma$. We will therefore refer to an equalizer prior as an *impartial* prior. The prior $\pi_0 \in \Pi_\infty$ is *minimax* if $\sup_{\tau \in \Phi} L(\tau, \pi_0) = \overline{W}$, where

$$\overline{W} = \inf_{\pi \in \Pi_\infty} \sup_{\tau \in \Phi} L(\tau, \pi)$$

is the upper value of the game. To obtain minimax solutions, we will adopt a standard game theory technique of searching for equalizer rules and showing that they are "extended Bayes" rules; see, for example, Chapter 5 of [4]. This is also the strategy used by Liang and Barron [19] for deriving minimax priors under the predictive regret (2.2) for location and scale families. In the present context the relevant result is given as Theorem 6.1 below.

Let $\Phi^+ \subset \Pi_\infty$ be the class of priors $\pi$ in $\Pi_\infty$ for which there exists a sequence $(\tau_k)$ of priors in $\Phi$ satisfying (i) $L(\tau_k, \pi) = \int L(\theta, \pi) \, d\tau_k(\theta)$ and (ii) $d(\tau_k, \pi) \to 0$. Since $L(\tau, \pi)$ is a proper scoring rule, each $\tau_k$ is a Bayes solution and, hence, $\Phi^+$ can be regarded as a class of extended Bayes solutions. If $\pi \in \Phi^+$ is an equalizer prior, then we can unambiguously define its predictive information as

$$\zeta(\pi) = \lim_{k \to \infty} \zeta(\tau_k)$$

for any sequence $\tau_k \in \Phi$ satisfying (i) and (ii) above. This is true since $L(\theta, \pi) = c$, say, for all $\theta \in \Theta$, and so for every such sequence we have $L(\tau_k, \pi) = c$ for all $k$ from (i). Therefore, from (4.10),

$$\zeta(\tau_k) = d(\tau_k, \pi) - c, \tag{6.2}$$

which tends to $-c$ as $k \to \infty$.

Finally, we define the class $U \subset \Pi_\infty$ of priors $\pi$ for which

$$\limsup_{n \to \infty} c_n \sup_{\theta \in \Theta} L_{Y|X}(\theta, \pi) < \infty \tag{6.3}$$

for every sequence $(m_n)$. Clearly, priors in $U^c$ have poor finite sample predictive behavior relative to Jeffreys' prior.

LEMMA 6.1. *Suppose that* $\pi \in C \cap U$, *that* R1, R2 *or* R3 *holds and that* $(m_n)$ *is any sequence satisfying the conditions in Theorem* 4.1(a), (b) *or* (c), *respectively. Then*

$$\sup_{\tau \in \Phi} L(\tau, \pi) \le \sup_{\theta \in \Theta} L(\theta, \pi).$$

PROOF. Let $\tau \in \Phi, \varepsilon > 0$ and choose a compact set $K \subset \Theta$ for which $\int_{K^c} d\tau(\theta) \le \varepsilon$. Then

$$L_{Y|X}(\tau, \pi) \le \sup_{\theta \in K} L_{Y|X}(\theta, \pi) + \varepsilon \sup_{\theta \in K^c} L_{Y|X}(\theta, \pi)$$

so that

$$L(\tau, \pi) = \limsup_{n \to \infty} c_n L_{Y|X}(\tau, \pi) \le \sup_{\theta \in K} L(\theta, \pi) + k\varepsilon$$

from (4.6) since $\pi \in C$, where $k = \limsup_{n \to \infty} c_n \sup_\theta L_{Y|X}(\theta, \pi) < \infty$ since $\pi \in U$. The result follows since $\varepsilon$ was arbitrary. $\square$

We now establish the following connection between equalizer and minimax priors.

THEOREM 6.1. *Suppose that* $\pi_0 \in \Phi^+ \cap C \cap U$ *is an equalizer prior, that* R1, R2 *or* R3 *holds with* $\pi = \pi_0$ *and that* $(m_n)$ *is any sequence satisfying the conditions in Theorem* 4.1(a), (b) *or* (c) *respectively. Then* $\pi_0$ *is minimax and* $\zeta(\pi_0) = \inf_{\tau \in \Phi} \zeta(\tau)$.

PROOF. Define

$$\underline{W} = \sup_{\tau \in \Phi} \inf_{\pi \in \Pi_\infty} L(\tau, \pi)$$

to be the lower value of the game. Then $\underline{W} \le \overline{W}$ is a standard result from game theory. Next, since $\pi_0$ is an equalizer prior, we have $L(\theta, \pi_0) = c$, say, for all $\theta \in \Theta$. Therefore, $\overline{W} = \inf_{\pi \in \Pi_\infty} \sup_{\tau \in \Phi} L(\tau, \pi) \leq \sup_{\tau \in \Phi} L(\tau, \pi_0) \leq \sup_{\theta \in \Theta} L(\theta, \pi_0) = c$ from Lemma 6.1 since $\pi_0 \in C \cap U$. Therefore, $\overline{W} \leq c$.

Since from Lemma 4.1 $L(\tau, \pi)$ is a $\Phi$-proper scoring rule, we have $\inf_{\pi \in \Pi_\infty} L(\tau, \pi) \geq L(\tau, \tau) = -\zeta(\tau)$ for every $\tau \in \Phi$. Therefore, $\underline{W} \geq -\inf_{\tau \in \Phi} \zeta(\tau)$. Since $\pi_0 \in \Phi^+$, there exists a sequence $(\tau_k)$ in $\Phi$ with $d(\tau_k, \pi_0) \to 0$. Therefore, since $\zeta(\tau_k) \geq \inf_{\tau \in \Phi} \zeta(\tau) \geq -\underline{W}$ and, from (6.2), $\zeta(\tau_k) \to -c$ as $k \to \infty$, we have $c \leq \underline{W}$. These relations give $\overline{W} \leq c \leq \underline{W}$ and it follows that $\underline{W} = c = \overline{W}$. The result now follows from the definitions of minimaxity and $\zeta(\pi_0)$. □

We see that, under the conditions of Theorem 6.1, the minimax prior $\pi_0$ has a natural interpretation of containing minimum predictive information about $Y$, since the infimum of the predictive information (4.9) is attained at $\tau = \pi_0$. Equivalently, $\pi_0$ is maximin since it maximizes the Bayes risk $-\zeta(\tau)$ of $\tau \in \Phi$ under (4.8) and, hence, is a *least favorable* prior under predictive loss. Notice also that Theorem 6.1 implies that $\sup_{\tau \in \Phi} L(\tau, \pi_0) = c$, regardless of the particular sequence $(m_n)$ used.

We note that for the assertion of Theorem 6.1 to hold we require that $\pi_0$ satisfies condition (6.3). There may exist a prior $\pi_1 \in U^c$ which appears to dominate the minimax prior $\pi_0$ on the basis of the asymptotic predictive loss function $L(\theta, \pi)$. However, this prior will possess poor penultimate asymptotic behavior since $L_{Y|X}(\theta, \pi)$ will be asymptotically unbounded. This will be reflected in the value of $\sup_{\tau \in \Phi} L(\tau, \pi)$, which will necessarily be greater than $\sup_{\theta \in \Theta} L(\theta, \pi)$. This phenomenon will be illustrated in Example 6.1.

COROLLARY 6.1. *Assume the conditions of Theorem 6.1 and additionally that $\pi_0$ is proper. Then if $\zeta(\pi_0) = -c$, where $c$ is the constant value of $L(\theta, \pi_0)$, then $\pi_0$ is minimax and $\zeta(\pi_0) = \inf_{\tau \in \Phi} \zeta(\tau)$.*

PROOF. Since $d(\pi_0, \pi_0) = 0$ and $\int L(\theta, \pi_0) \, d\pi_0(\theta) = c = -\zeta(\pi_0) = L(\pi_0, \pi_0)$, it follows on taking $\tau_k = \pi_0$ that $\pi_0 \in \Phi^+$. The result now follows from Theorem 6.1. □

Suppose that $\pi_0 \in C \cap U$ is an improper equalizer prior. One way to show that $\pi_0 \in \Phi^+$ is to construct a sequence $(\tau_k)$ of priors in $\Gamma$ for which $d(\tau_k, \pi_0) \to 0$, where $d(\tau, \pi_0)$ is given by formula (5.14). As noted just prior to Lemma 5.5, the condition $L(\tau_k, \pi_0) = \int L(\theta, \pi_0) \, d\tau_k(\theta)$ is automatically satisfied when $\tau_k \in \Gamma$.

We consider first the case $p = 1$. In this case it turns out that Jeffreys' prior is a minimax solution, and, hence, the assertion at the end of Example 5.1. Let $\mathcal{H}$ be the class of probability density functions $h$ on $(-1, 1)$ possessing second-order continuous derivatives and that satisfy $h(-1) = h'(-1) =$

$h''(-1) = h(1) = h'(1) = h''(1) = 0$ and

(6.4) $$\int_{-1}^{1} \{g'(u)\}^2 h(u)\, du < \infty,$$

where $g(u) = \log h(u)$; that is, the Fisher information associated with $h$ is finite. The class $\mathcal{H}$ is nonempty, since the density of the random variable $U = 2V - 1$, where $V$ is any beta $(a, b)$ density with $a, b > 3$, satisfies these conditions.

COROLLARY 6.2. *Suppose that $p = 1$. Then Jeffreys' prior is minimax and $\zeta(\pi^J) = \inf_{\tau \in \Phi} \zeta(\tau)$.*

PROOF. Since $L(\theta, \pi^J) = 0$, Jeffreys' prior is an equalizer prior. We therefore need to show that $\pi^J \in \Phi^+ \cap C \cap U$. Recall that $\pi^J \in C$ was an assumption made in Section 4. Also, since $L_{Y|X}(\theta, \pi^J) = 0$ for all $n$ from (2.5), $\pi^J \in U$.

If $\pi^J$ is proper, the result now follows immediately from Corollary 6.1 since $\zeta_{Y|X}(\pi^J) = 0$ for all $n$. Suppose then that $\pi^J$ is improper. Without loss of generality, we assume that $i(\theta) = 1$, so that Jeffreys' prior is uniform. Since $\pi^J$ is improper, without loss of generality we take $\Theta$ to be either $(-\infty, \infty)$ or $(0, \infty)$ by a suitable linear transformation. Now let $U$ be a random variable with density $h \in \mathcal{H}$.

Suppose first that $\Theta = (-\infty, \infty)$ and let $\tau_k$ be the density of $\theta = kU$. Clearly, $\tau_k \in \Gamma$, $\tau_k$ has support $[-k, k]$ and $\mu'_k(\theta) = g'(u)/k$, where $\mu_k = \log \tau_k$ and $u = \theta/k$. Therefore, from (5.14),

$$d(\tau_k, \pi^J) = \frac{1}{k^2} E\{g'(U)\}^2 \to 0$$

as $k \to \infty$ from (6.4) so that $\pi^J \in \Phi^+$. The result now follows from Theorem 6.1.

Next suppose that $\Theta = (0, \infty)$ and let $\tau_k$ be the density of $\theta = k(U+1)+1$. Then $\tau_k \in \Gamma$, $\tau_k$ has support $[1, 2k+1]$ and $\mu'_k(\theta) = g'(u)/k$, where $u = (\theta - 1)/k - 1$. Therefore, from (5.14),

$$d(\tau_k, \pi^J) = \frac{1}{k^2} E\{g'(U)\}^2 \to 0$$

as $k \to \infty$ from (6.4), so that $\pi^J \in \Phi^+$ and again the result follows from Theorem 6.1. □

EXAMPLE 6.1. Bernoulli model. Here Jeffreys' prior is the beta $(1/2, 1/2)$ distribution, which is therefore minimax from Corollary 6.2. The underlying Bernoulli probability mass function is $f(x|\theta) = \theta^x (1-\theta)^{1-x}, x = 0, 1, 0 < \theta <$

1. Let $\pi^a$ be the density of the beta $(a,a)$ distribution, where $a > 0$. It is straightforward to check from (5.12) that

$$L(\theta, \pi^a) = \left(a - \frac{1}{2}\right)\left\{-4\left(a - \frac{1}{2}\right) + \frac{a - 3/2}{\theta(1-\theta)}\right\},$$

from which we see that $L(\theta, \pi_1) = -4$, where $\pi_1 = \pi^{3/2}$, the beta $(\frac{3}{2}, \frac{3}{2})$ distribution. Hence, the prior $\pi_1$ would appear to dominate Jeffreys' prior. In view of Corollary 6.2, however, we conclude that condition (6.3) must break down for this prior. Indeed, it can be shown directly that $c_n L_{Y|X}(0, \pi_1)$ is an increasing function of $m$ for fixed $n$ and that, when $m = 1$, we have $c_n L_{Y|X}(0, \pi_1) = n + O(1)$. By the continuity of $L_{Y|X}(\theta, \pi_1)$ in $(0, 1)$, it follows that $c_n \sup_\theta L_{Y|X}(\theta, \pi_1) \to \infty$ as $n \to \infty$ for every sequence $(m_n)$ and so $\pi_1 \notin U$. Therefore, $\pi_1$ exhibits poor finite sample predictive behavior relative to Jeffreys' prior for values of $\theta$ close to 0 or 1.

It is of some interest to compare this behavior with the asymptotic minimax analysis under the prior predictive regret (4.1). Under (4.1), Jeffreys' prior is asymptotically maximin [8], but not minimax due to its poor boundary risk behavior. However, a sequence of priors converging to Jeffreys' prior can be constructed that is asymptotically minimax [26]. Under our posterior predictive regret criterion, Jeffreys' prior is both maximin and minimax. In particular, it follows that it is not possible to modify the beta $(\frac{3}{2}, \frac{3}{2})$ distribution at the boundaries to make it asymptotically minimax.

In the examples below our strategy for identifying a minimax prior will be to consider a suitable class of candidate priors in $C$, compute the predictive loss (5.12), identify the subclass of equalizer priors in $U$ and choose the prior $\pi_0$ in this subclass, assuming it is nonempty, with minimum constant loss. Clearly, $\pi_0$ will be minimax over this subclass of equalizer priors. If, in addition, it can be shown that $\pi_0 \in \Phi^+$, then the conditions of Theorem 6.1 hold and $\pi_0$ is minimax over $\Phi$. In particular, we will see that in dimensions greater than one, although Jeffreys' prior is necessarily impartial, it may not be minimax. This is not surprising, since we know that in the special case of transformation models the right Haar measure is the best invariant prior under posterior predictive loss (see Section 2). Exact minimax solutions for Examples 6.2 and 6.3 under the predictive regret (2.2) have recently been obtained by Liang and Barron [19]. Finally, all these examples are sufficiently regular for the strongest form R3 of remainder to hold for the priors $\pi_0$ that are obtained. Hence, from Theorem 4.1(c), all the results will apply for an arbitrary amount of prediction.

EXAMPLE 6.2. Normal model with unknown mean and variance. Here $X \sim N(\beta, \sigma^2)$ and $\theta = (\beta, \sigma)$. We will show that the prior $\pi_0(\theta) \propto \sigma^{-1}$ is minimax. This is Jeffreys' independence prior, or the right Haar measure under the group of affine transformations of the data.

Consider the class of improper priors $\pi^a(\theta) \propto \sigma^{-a}$ on $\Theta$, where $a \in \mathcal{R}$. Transforming to $\phi = (\beta, \lambda)$, where $\lambda = \log \sigma$, these priors become $\pi^a(\phi) \propto \exp\{-(a-1)\lambda\}$ in the $\phi$-parameterization. Here we find that $i(\phi) = \mathrm{diag}(e^{-2\lambda}, 2)$. Since $\rho^a(\phi) = \log \pi^a(\phi) = -(a-1)\lambda$, it follows immediately from (5.2) that $A(\phi, \pi^a) = \frac{1}{2}(a-1)^2$. Furthermore, since $|i(\phi)| = 2e^{-2\lambda}$, we have $\pi^J(\phi) \propto e^{-\lambda} = \pi^2(\phi)$ so that $A(\phi, \pi^J) = \frac{1}{2}$. It now follows from (5.12) that $L(\phi, \pi^a) = \frac{1}{2}\{(a-1)^2 - 1\}$. Therefore, all priors in this class are equalizer priors and $L(\phi, \pi^a)$ attains its minimum value in this class when $a = 1$, which corresponds to $\pi_0(\phi) \propto 1$, or $\pi_0(\theta) \propto \sigma^{-1}$ in the $\theta$-parameterization. Note that the minimum value $-\frac{1}{2} < 0$, which is the loss under Jeffreys' prior.

We now show that $\pi_0 \in \Phi^+ \cap C \cap U$. Clearly, $\pi_0 \in C$, while $\pi_0 \in U$ follows because $L_{Y|X}(\theta, \pi_0)$ is constant for all $n$ since $\pi_0$ is invariant under the transitive group of transformations of $\Theta$ induced by the group of affine transformations of the observations (see Section 2). It remains to show that $\pi_0 \in \Phi^+$. Let $U_1, U_2$ be independent random variables with common density $h \in \mathcal{H}$ and let $\tau_k$ be the joint density of $\phi = (\beta, \lambda)$, where $\beta = k_1 U_1, \lambda = k_2 U_2$ and $k_1, k_2$ are functions of $k$ to be determined. Let $\mu_k = \log \tau_k$. Then $\mu_{kr} = k_r^{-1} g'(U_r), r = 1, 2$, where $g = \log h$. Write $\alpha = \int_{-1}^{1} \{g'(u)\}^2 h(u)\, du < \infty$ since $h \in \mathcal{H}$. Since $\rho_0(\phi) = \log \pi_0(\phi)$ is constant, it follows from (5.14) that

$$d(\tau_k, \pi_0) = E[k_1^{-2} e^{2\lambda}\{g'(U_1)\}^2 + \tfrac{1}{2} k_2^{-2}\{g'(U_2)\}^2] \leq \alpha\{k_1^{-2} e^{2k_2} + \tfrac{1}{2} k_2^{-2}\},$$

since $\lambda \leq k_2$. Now take $k_1 = ke^k, k_2 = k$. Then $d(\tau_k, \pi_0) \leq \frac{3\alpha}{2k^2} \to 0$ as $k \to \infty$ and, hence, $\pi_0 \in \Phi^+$. It now follows from Theorem 6.1 that $\pi_0$ is minimax and that $\zeta(\pi_0) = \frac{1}{2}$.

EXAMPLE 6.3. Normal linear regression. Here $X_i \sim N(z_i^T \beta, \sigma^2)$, $i = 1, \ldots, n$, where $Z_n = (z_1, \ldots, z_n)^T$ is an $n \times q$ matrix of rank $q \geq 1$ and $\theta = (\beta, \sigma)$. Using a similar argument to that in Example 6.2, we can show that again Jeffreys' independence prior, or the right Haar measure, $\pi_0(\theta) \propto \sigma^{-1}$ is minimax.

Since the variables are not identically distributed in this example, it is not covered by the asymptotic theory of Sections 4 and 5. However, under suitable stability assumptions on the sequence $(z_i)$ of regressor variables, at least that $V_n \equiv n^{-1} Z_n^T Z_n$ is uniformly bounded away from zero and infinity, then a version of Theorem 5.1 will apply.

Proceeding as in Example 6.2, we again consider the class of priors $\pi^a(\theta) \propto \sigma^{-a}$ on $\Theta$, where $a \in \mathcal{R}$. Transforming to $\phi = (\beta, \lambda)$, where $\lambda = \log \sigma$, these priors become $\pi^a(\phi) \propto \exp\{-(a-1)\lambda\}$. Here we find that $i_n(\phi) = \mathrm{diag}(e^{-2\lambda} V_n, 2)$ and, exactly as in Example 6.2, we obtain $A(\phi, \pi^a) = \frac{1}{2}(a-1)^2$. Here $|i_n(\phi)| = 2|V_n| e^{-2q\lambda}$ so $\pi^J(\phi) \propto e^{-q\lambda} = \pi^{q+1}(\phi)$ for all $n$, giving $A(\phi, \pi^J) = \frac{1}{2} q^2$ and, hence, $L(\phi, \pi^a) = \frac{1}{2}\{(a-1)^2 - q^2\}$. Therefore, all priors in this class are equalizer priors and $L$ attains its minimum value in this class when $a = 1$, which corresponds to $\pi_0(\phi) \propto 1$, or $\pi_0(\theta) \propto \sigma^{-1}$ in the $\theta$-parameterization. Notice that the drop in predictive loss increases as the square of the number $q$ of regressors in the model. Note also that the ratio $|i_n|^{-1}|i_{n+1}|$ is free from $\theta$, so that a version of Theorem 4.1 will hold.

Exactly as in Example 6.2, $\pi_0 \in C \cap U$ and it remains to show that $\pi_0 \in \Phi^+$. Let $p = q+1$ and $U_j, j = 1, \ldots, p$, be independent random variables with common density $h \in \mathcal{H}$. With the same definitions as in Example 6.2, let $\beta_r = k_1 U_r, r = 1, \ldots, q, \lambda = k_2 U_p$, so that $\mu_{kr} = k_1^{-1} g'(U_r), r = 1, \ldots, q, \mu_{kp} = k_2^{-1} g'(U_p)$. Then it follows from (5.14) that, with the summations over $r$ and $s$ running from 1 to $q$,

$$\begin{aligned} d(\tau_k, \pi_0) &= E\{e^{2\lambda} V_n^{rs} \mu_{kr} \mu_{ks} + \tfrac{1}{2} \mu_{kp}^2\} \\ &= E\{k_1^{-2} e^{2\lambda} V_n^{rs} g'(U_r) g'(U_s) + \tfrac{1}{2} k_2^{-2} g'(U_p)^2\} \\ &\leq \alpha\{k_1^{-2} e^{2k_2} \mathrm{trace}(V_n^{-1}) + \tfrac{1}{2} k_2^{-2}\}, \end{aligned}$$

using $\int_{-1}^{1} g'(u) h(u)\, du = 0$. Now take $k_1 = ke^k, k_2 = k$. Then, as before, $d(\tau_k, \pi_0) \to 0$ as $k \to \infty$ and, hence, $\pi_0 \in \Phi^+$. It follows from Theorem 6.1 that $\pi_0$ is minimax and $\zeta(\pi_0) = \frac{q^2}{2}$.

Interestingly, we note that the priors $\pi_0$ identified in Examples 6.2 and 6.3 also give rise to minimum predictive coverage probability bias; see [12]. The next example is more challenging and illustrates the difficulties associated with finding minimax priors more generally.

EXAMPLE 6.4. Multivariate normal. Here $X \sim N_q(\mu, \Sigma)$, with $\theta$ comprising all elements of $\mu$ and $\Sigma$. Write $\Sigma^{-1} = T'T$, where $T = (t_{ij})$ is a lower triangular matrix satisfying $t_{ii} > 0$. Let $\mu = (\mu_1, \ldots, \mu_q)', \psi_i = t_{ii}, 1 \leq i \leq q, \psi = (\psi_1, \ldots, \psi_q)', \beta_{ij} = t_{ii}^{-1} t_{ij}, 1 \leq j < i \leq q$ and $\beta^{(i)} = (\beta_{i1}, \ldots, \beta_{i,i-1})', 2 \leq i \leq q$. Then $\gamma = (\psi', \beta^{(2)'}, \ldots, \beta^{(q)'}, \mu')'$ is a one-to-one transformation of $\theta$. The loglikelihood is

$$l(\gamma) = \sum_{i=1}^{q} \log \psi_i - \tfrac{1}{2}\left[\sum_{i=1}^{q} \psi_i^2 \left\{\sum_{j=1}^{i} \beta_{ij}(x_j - \mu_j)\right\}^2\right],$$

writing $\beta_{ii} = 1, i = 1, \ldots, q$. One then finds that the information matrix $i(\gamma)$ is block diagonal in $\psi_1, \ldots, \psi_q, \beta^{(2)'}, \ldots, \beta^{(q)'}, \mu'$ and is given by

$$\mathrm{diag}(2\psi_1^{-2}, \ldots, 2\psi_q^{-2}, \psi_2^2 \Sigma_{11}, \ldots, \psi_q^2 \Sigma_{q-1,q-1}, \Sigma^{-1}),$$

where $\Sigma_{ii}$ is the submatrix of $\Sigma$ corresponding to the first $i$ components of $X$. Using the fact that $|\Sigma_{ii}| = \prod_{j=1}^{i} \psi_j^{-2}, i = 1, \ldots, q$, we obtain $|i(\gamma)| = 2^q \prod_{i=1}^{q} \psi_i^{4i-2q-2}$.

Consider the class of priors $\pi^a(\theta) \propto |\Sigma|^{-(q+2-a)/2}$ on $\Theta$, where $a \in \mathcal{R}$. In the $\gamma$-parameterization, this class becomes $\pi^a(\gamma) \propto \prod_{i=1}^{q} \psi_i^{2i-q-a-1}$. Noting that the case $a = 0$ is Jeffreys' prior, it is straightforward to show from (5.12) that $L(\gamma, \pi) = \frac{q}{2}\{(a-1)^2 - 1\}$. Therefore, all priors in this class are equalizer priors and $L$ attains its minimum value within this class when $a = 1$. From invariance considerations via affine transformations of $X$, it can be shown that these priors are also equalizer priors for finite $n$ and, hence, are all in the class $U$. These results therefore suggest that the right Haar prior $\pi_0(\theta) \propto |\Sigma|^{-(q+1)/2}$ arising from the affine group is minimax. However, in this example it does not appear to be possible to approximate $\pi_0$ by a sequence of compact priors, as was done in the previous examples. We conjecture, however, that $\pi_0$ can be approximated by a suitable sequence of proper priors so that Theorem 6.1 will give the minimaxity of $\pi_0$, but we have been unable to demonstrate this. This example does show, however, that Jeffreys' prior is dominated by $\pi_0$.

Interestingly, further analysis reveals that the prior $\pi_1(\gamma) \propto \prod_{i=1}^{q} \psi_i^{-1}$ is also an equalizer prior and that it dominates $\pi_0$. In the $\theta$-parameterization this prior becomes $\pi_1(\theta) \propto \{\prod_{i=1}^{q} |\Sigma_{ii}|\}^{-1}$. However, this prior is seen to be noninvariant under nonsingular transformation of $X$ and, furthermore, does not satisfy the boundedness condition (6.3).

In the case $q = 2$, in the parameterization $\phi = (\mu_1, \mu_2, \sigma_1, \sigma_2, \rho)$, where $\sigma_i$ is the standard deviation of $X_i, i = 1, 2$, and $\rho = \text{Corr}(X_1, X_2)$, Jeffreys' prior and $\pi_0$ become, respectively,

$$\pi^J(\phi) \propto \sigma_1^{-2}\sigma_2^{-2}(1-\rho^2)^{-2},$$
$$\pi_0(\phi) \propto \sigma_1^{-1}\sigma_2^{-1}(1-\rho^2)^{-3/2}.$$

Therefore (see the paragraph below), $\pi_0$ is Jeffreys' "two-step" prior. In the context of our predictive set-up, marginalization issues correspond to predicting only certain functions of the future data $Y = (X_{n+1}, \ldots, X_{n+m})$. In general, the associated minimax predictive prior will differ from that for the problem of predicting the entire future data $Y$ unless the selected statistics just form a sufficiency reduction of $Y$. Such questions will be explored in future work. Thus, if we were only interested in predicting the correlation coefficient of a future set of bivariate data, then we might start with the observed correlation as the data $X$ and use Jeffreys' prior in this single parameter case, which is $\pi(\rho) \propto (1-\rho^2)^{-1}$. For further discussion and references on the choice of prior in this example, see [6], page 363.

Finally, we note the corresponding result for general $q$ in the case $\mu$ known. Again, considering the class of priors $\pi^a(\theta) \propto |\Sigma|^{-(q+2-a)/2}$ on $\Theta$, we find that the optimal choice is $a = 1$, so $\pi_0$ is as given above and in this case coincides with Jeffreys' prior. This was also shown to be a predictive probability matching prior in [12] in the case $q = 2$.

Under the conditions of Theorem 6.1, it is possible to change the base measure from Jeffreys' prior to $\pi_0$, since $\pi_0$ is neutral with respect to $\pi^J$ under $L(\theta, \pi)$. Denoting quantities with respect to the base measure $\pi_0$ with a zero subscript, since $L(\theta, \pi_0) = c \leq 0$ and $\zeta(\pi_0) = -c$, we have, for $\pi \in \Pi_\infty$,

$$L_0(\theta, \pi) = L(\theta, \pi) - L(\theta, \pi_0) = L(\theta, \pi) - c$$

and for $\tau \in \Phi$,

$$\zeta_0(\tau) = \zeta(\tau) + c.$$

Therefore, with respect to the base measure $\pi_0$, the predictive loss under $\pi_0$ becomes $L_0(\theta, \pi_0) = 0$ and the minimum predictive information, attained at $\pi = \pi_0$, is zero.

**7. Discussion.** In this paper we have obtained an asymptotic predictive loss function that reflects the finite sample size predictive behavior of alternative priors when the sample size is large for arbitrary amounts of prediction. This loss function is related to that in [14] for the comparison of estimative predictive distributions based on Bayes estimators. It can be used to derive nonsubjective priors that are impartial, minimax and maximin, which is equivalent here to minimizing a measure of the predictive information contained in a prior. In dimensions greater than one, unlike an analysis based on prior predictive regret, the maximin prior may not be Jeffreys' prior. A number of examples have been given to illustrate these ideas.

As discussed in [23], as model complexity increases, it becomes more difficult to make sensible prior assignments, while at the same time the effect of the prior specification on the final inference of interest becomes more pronounced. It is therefore important to have sound methodology available for the construction and implementation of priors in the multiparameter case. We believe that our preliminary analysis of the posterior predictive regret (2.1) indicates that it should be a valuable tool for such an enterprise. More extensive analysis is now required, particularly aimed at developing general methods of finding exact and approximate solutions for the practical implementation of this work and investigating connections with predictive coverage probability bias. Local priors (see, e.g., [23, 24]) are expected to play a role. It would also be interesting to develop asymptotically impartial minimax posterior predictive loss priors for dependent observations and for various classes of nonregular problems. In particular, all the definitions in Section 2 for nonasymptotic settings will apply and could be used to explore predictive behavior numerically.

**Acknowledgments.** We would like to thank two referees and an Associate Editor for their constructive comments and suggestions for improving the clarity of this paper.

## REFERENCES

[1] AITCHISON, J. (1975). Goodness of prediction fit. *Biometrika* **62** 547–554. MR0391353
[2] AKAIKE, H. (1978). A new look at the Bayes procedure. *Biometrika* **65** 53–59. MR0501450
[3] BARRON, A. R. (1999). Information-theoretic characterization of Bayes performance and the choice of priors in parametric and nonparametric problems. In *Bayesian Statistics 6* (J. M. Bernardo, J. O. Berger, A. P. Dawid and A. F. M. Smith, eds.) 27–52. Oxford Univ. Press, New York. MR1723492
[4] BERGER, J. O. (1985). *Statistical Decision Theory and Bayesian Analysis*, 2nd ed. Springer, New York. MR0804611
[5] BERNARDO, J. M. (1979). Reference posterior distributions for Bayesian inference (with discussion). *J. Roy. Statist. Soc. Ser. B* **41** 113–147. MR0547240
[6] BERNARDO, J. M. and SMITH, A. F. M. (1994). *Bayesian Theory*. Wiley, Chichester. MR1274699
[7] CLARKE, B. S. and BARRON, A. R. (1990). Information-theoretic asymptotics of Bayes methods. *IEEE Trans. Inform. Theory* **36** 453–471. MR1053841
[8] CLARKE, B. S. and BARRON, A. R. (1994). Jeffreys' prior is asymptotically least favorable under entropy risk. *J. Statist. Plann. Inference* **41** 37–60. MR1292146
[9] CLARKE, B. and YUAN, A. (2004). Partial information reference priors: Derivation and interpretations. *J. Statist. Plann. Inference* **123** 313–345. MR2062985
[10] COVER, T. M. and THOMAS, J. A. (1991). *Elements of Information Theory*. Wiley, New York. MR1122806
[11] DATTA, G. S. and MUKERJEE, R. (2004). *Probability Matching Priors*: *Higher Order Asymptotics. Lecture Notes in Statist.* **178**. Springer, New York. MR2053794
[12] DATTA, G. S., MUKERJEE, R., GHOSH, M. and SWEETING, T. J. (2000). Bayesian prediction with approximate frequentist validity. *Ann. Statist.* **28** 1414–1426. MR1805790
[13] GHOSH, J. K. and MUKERJEE, R. (1991). Characterization of priors under which Bayesian and frequentist Bartlett corrections are equivalent in the multiparameter case. *J. Multivariate Anal.* **38** 385–393. MR1131727
[14] HARTIGAN, J. A. (1998). The maximum likelihood prior. *Ann. Statist.* **26** 2083–2103. MR1700222
[15] KASS, R. E., TIERNEY, L. and KADANE, J. (1990). The validity of posterior expansions based on Laplace's method. In *Bayesian and Likelihood Methods in Statistics and Econometrics* (S. Geisser, J. S. Hodges, S. J. Press and A. Zellner, eds.) 473–488. North-Holland, Amsterdam.
[16] KOMAKI, F. (1996). On asymptotic properties of predictive distributions. *Biometrika* **83** 299–313. MR1439785
[17] KUBOKI, H. (1998). Reference priors for prediction. *J. Statist. Plann. Inference* **69** 295–317. MR1631332
[18] LAUD, P. W. and IBRAHIM, J. G. (1995). Predictive model selection. *J. Roy. Statist. Soc. Ser. B* **57** 247–262. MR1325389
[19] LIANG, F. and BARRON, A. R. (2004). Exact minimax strategies for predictive density estimation, data compression, and model selection. *IEEE Trans. Inform. Theory* **50** 2708–2726. MR2096988
[20] MCCULLAGH, P. (1987). *Tensor Methods in Statistics*. Chapman and Hall, London. MR0907286
[21] SUN, D. and BERGER, J. O. (1998). Reference priors with partial information. *Biometrika* **85** 55–71. MR1627242

[22] SWEETING, T. J. (1996). Approximate Bayesian computation based on signed roots of log-density ratios (with discussion). In *Bayesian Statistics 5* (J. M. Bernardo, J. O. Berger, A. P. Dawid and A. F. M. Smith, eds.) 427–444. Oxford Univ. Press, New York. MR1425418
[23] SWEETING, T. J. (2001). Coverage probability bias, objective Bayes and the likelihood principle. *Biometrika* **88** 657–675. MR1859400
[24] SWEETING, T. J. (2005). On the implementation of local probability matching priors for interest parameters. *Biometrika* **92** 47–57. MR2158609
[25] TIBSHIRANI, R. (1989). Noninformative priors for one parameter of many. *Biometrika* **76** 604–608. MR1040654
[26] XIE, Q. and BARRON, A. R. (1997). Minimax redundancy for the class of memoryless sources. *IEEE Trans. Inform. Theory* **43** 646–657.

G. S. DATTA  
DEPARTMENT OF STATISTICS  
UNIVERSITY OF GEORGIA  
ATHENS, GEORGIA 30602-1952  
USA  
E-MAIL: gauri@stat.uga.edu

M. GHOSH  
DEPARTMENT OF STATISTICS  
UNIVERSITY OF FLORIDA  
GAINESVILLE, FLORIDA 32611-8545  
USA  
E-MAIL: ghoshm@stat.ufl.edu

T. J. SWEETING  
DEPARTMENT OF STATISTICAL SCIENCE  
UNIVERSITY COLLEGE LONDON  
GOWER STREET  
LONDON, WC1E 6BT  
UNITED KINGDOM  
E-MAIL: trevor@stats.ucl.ac.uk